\documentstyle[11pt]{article}
\title{Time-Frequency Localization and the Gabor Transform}
\author{Jayakumar Ramanathan\thanks{Supported by a Faculty Research
and Creative Projects Fellowship Award from Eastern Michigan
University.}\\
Department of Mathematics\\
Eastern Michigan University\\
Ypsilanti MI 48197
\and
Pankaj Topiwala\thanks{Supported by the MITRE Sponsored Research
program.  P. T. would like to thank H. Dym, A. J. E. M. Janssen,
and P. Flandrin for valuable discussions.}\\
The MITRE Corporation\\
Burlington Road\\
Bedford MA 01730}

\newtheorem{lemma}{Lemma}
\newtheorem{proposition}{Proposition}
\newtheorem{theorem}{Theorem}

\begin{document}

\maketitle

\begin{abstract} The time-frequency content of a signal can be measured by the
Gabor transform or windowed Fourier transform. This is a function defined on
phase space that is computed by taking the Fourier transform of the product
of the signal against a translate of a fixed window. The problem of finding
signals with Gabor transform that are maximally concentrated within a
given region of phase space is discussed. It is well known that such problems
give rise to an eigenvalue problem for an associated self-adjoint, positive
{\em concentration} operator that has its spectrum contained in the unit
interval. In this paper, the asymptotic behavior of these eigenvalues as the
concentration region gets large is studied. The smoothness of the
eigenfunctions is also examined.
\end{abstract}

\section{Introduction}

A time-frequency shift of a square integrable function
$f(t)\in L^2({\bf R})$ is defined by:
\[
\rho(\tau,\sigma)f(t) = \exp(\pi i \tau\sigma) \exp(2\pi i \sigma t) f(t +
\tau).
\]
The first factor is not essential and is present in order to simplify
many of the susequent formulae. Ignoring this, the formula clearly has the
structure of a time shift of $\tau$ units together with a frequency shift
of $\sigma$ units. The Gabor transform or windowed Fourier transform depends
on the choice of a window function $\phi(t)\in L^2({\bf R})$ and is
defined by
\begin{eqnarray*}
{\cal S}_\phi f(\tau,\sigma) &=& \langle f, \rho(\tau,\sigma)\phi\rangle\\
&=& \int f(t) \overline{\rho(\tau,\sigma)\phi(t)}\,dt\\
\end{eqnarray*}
for any $f\in L^2({\bf R})$.
Here $\langle\bullet,\bullet\rangle$ denotes the standard inner product
on $ L^2({\bf R})$. The Fourier transform of $f(t)$ is defined by
\[
\hat f(\sigma) = {\cal F}f(\sigma) = \int f(t) e^{-2\pi i t\sigma}\,dt.
\]
One easily checks that
\begin{equation}
{\cal F}\rho(\tau,\sigma)f = \rho(-\sigma,\tau)\hat f.
\label{eqn:Fshift}
\end{equation}

The following proposition collects some formulae that  will be useful in
the sequel.

\begin{proposition}
Let $f,\phi\in L^2({\bf R})$.
\begin{enumerate}
\item $\int\int |{\cal S}_\phi f(\tau,\sigma)|^2\,d\tau d\sigma = \|f\|^2
\|\phi\|^2$.
\item ${\cal S}_\phi f(\tau,\sigma) = {\cal S}_{\hat \phi} \hat
f(-\sigma,\tau)$ where $\hat f$ and $\hat \phi$ denote the Fourier
transform.
\end{enumerate}
\label{spectrogram_identities}
\end{proposition}
{\sc Proof:} The first item is a well known special case of Moyal's
identity. See \cite{Folland}. The second item is proved by using
Plancherel's theorem together with equation~\ref{eqn:Fshift}:
\begin{eqnarray*}
{\cal S}_\phi f(\tau,\sigma) &=& \langle
f,\rho(\tau,\sigma)\phi\rangle\\
&=& \langle {\cal F}f,{\cal F}\rho(\tau,\sigma)\phi\rangle\\
&=& {\cal S}_{\hat \phi} \hat f(-\sigma,\tau)
\end{eqnarray*}
$\Box$

Let $\Omega \subset {\bf R}^2$ be a region within the
time-frequency plane. For the remainder of the exposition, $\phi$ will be
assumed to be normalized to have $L^2$ norm equal to one: $\|\phi\|^2 = 1$.
We will measure the energy of $f$ contained within $\Omega$ by
\[
{\cal E}_\Omega(f) = \int\int_\Omega
|{\cal S}_\phi f(\tau,\sigma)|^2\,d\tau d\sigma.
\]
Associated to the energy is the Hermitian symmetric form defined by:
\[
{\cal E}_\Omega(f,g) = \int\int_\Omega \langle
f,\rho(\tau,\sigma)\phi\rangle
\langle\rho(\tau,\sigma)\phi,g\rangle\,d\tau d\sigma
\]
for any $f,g\in L^2({\bf R})$. Proposition~\ref{spectrogram_identities}
insures that ${\cal E}$ is bounded in the following sense:
\begin{equation}
{\cal E}_\Omega(f,g) \le \|f\|\ \|g\|.
\label{eqn:bound}
\end{equation}
A standard duality argument produces a bounded, positive, self-adjoint operator
${\cal C}_\Omega:L^2({\bf R})\to L^2({\bf R})$ such that
\[
{\cal E}_\Omega(f,g) = \langle {\cal C}_\Omega f,g\rangle
\]
for all $f,g \in L^2({\bf R})$. Let $\psi_1,\psi_2,\cdots$ be a complete
orthonormal
basis of $L^2({\bf R})$  and $\lambda_1,\lambda_2,\cdots$ be a decreasing
sequence of nonnegative real numbers for which
\[
{\cal C}_\Omega\psi_k = \lambda_k \psi_k.
\]
The bound in equation~\ref{eqn:bound} forces
$1\ge\lambda_1\ge\lambda_2\ge\cdots$.

It is clear that the unit norm signal with the maximum energy within
$\Omega$ is $\psi_1$. Moreover, this maximum energy is given by
$\lambda_1$. More generally, the quantity
\[
n(\lambda) = {\rm card}\{k: \lambda_k > \lambda\}
\]
determines the maximum dimension of a subspace $V \in L^2({\bf R})$ of signals
for which
\[
{\cal E}_\Omega(f) \ge \lambda \|f\|^2
\]
for all $f\in V$. Thus the distribution of these eigenvalues as well as
the smoothness of the corresponding eigenvalues are of
interest within the context of signal processing.

\section{A Reproducing Kernel Hilbert Space}

It will be convenient to analyse the distribution of the eigenvalues
using the reproducing kernel Hilbert space setting described in this
section.

Proposition~\ref{spectrogram_identities} implies that ${\cal S}_\phi$ is
an isometry of $L^2({\bf R})$ onto  a closed subspace ${\cal V}_\phi \subset
L^2({\bf R}^2)$. Let ${\cal P}_\phi$ be the orthogonal projection from
$L^2({\bf R}^2)$ to ${\cal V}_\phi$. Let $F \in {\cal V}_\phi$ and $f =
{\cal S}_\phi^{-1}F$. Then,
applying the Cauchy-Schwarz inequality, one has
\[
|F(\tau_0,\sigma_0)| = \langle f,\rho(\tau_0,\sigma_0)\phi\rangle
\le  \|f\|\ \|\rho(\tau_0,\sigma_0)\phi\| = \|F\|.
\]
Therefore point evaluation at a given $(\tau_0,\sigma_0)$ of functions in
${\cal V}_\phi$ is a bounded linear functional.
Let $K_{(\tau_0,\sigma_0)}(\tau,\sigma) =
K(\tau_0,\sigma_0;\tau,\sigma)$ denote that element of
${\cal V}_\phi$ for which
\[
F(\tau_0,\sigma_0) = \int\int
F(\tau,\sigma) \overline{K_{(\tau_0,\sigma_0)}(\tau,\sigma)}\,d\tau d\sigma.
\]
In order to identify the kernel explicitly one only has to notice that
\[
F(\tau_0,\sigma_0) = \langle{\cal
S}_\phi^{-1}F,\rho(\tau_0,\sigma_0)\phi\rangle.
\]
Since ${\cal S}_\phi$ is an isometry,
\begin{eqnarray*}
F(\tau_0,\sigma_0) &=& \langle F,{\cal
S}_\phi\rho(\tau_0,\sigma_0)\phi\rangle\\
&=& \int\int F(\tau,\sigma)
\overline{\langle
\rho(\tau_0,\sigma_0)\phi,\rho(\tau,\sigma)\phi\rangle}\,d\tau d\sigma\\
\end{eqnarray*}
Consequently,
\[
K(\tau_0,\sigma_0;\tau,\sigma) = \langle \rho(\tau_0,\sigma_0)\phi,
\rho(\tau,\sigma)\phi\rangle.
\]
Finally, if $F \in {\cal V}^\perp_\phi$ we have
\[
\int\int
F(\tau,\sigma)\overline{K_{(\tau_0,\sigma_0)}(\tau,\sigma)}\,d\tau
d\sigma = 0
\]
since $K_{(\tau_0,\sigma_0)}(\tau,\sigma) \in {\cal V}_\phi$.
These remarks are summarized in the following proposition.

\begin{proposition}
For any $F \in L^2({\bf R}^2)$ the projection ${\cal P}_\phi F$ is computed by
\[
{\cal P}_\phi F(\tau,\sigma) = \int\int
F(\tau',\sigma') \overline{K_{(\tau,\sigma)}(\tau',\sigma')}\,d\tau'
d\sigma'.
\]
\label{proj_formula}
\end{proposition}

Our next task is to use the isometry ${\cal S}_\phi$ to identify the
concentration operator ${\cal C}_\Omega$ explicitly. In particular, the
next proposition provides a formula for ${\cal S}_\phi{\cal C}_\Omega
{\cal S}_\phi^{-1}$ on ${\cal V}_\phi$.

\begin{proposition}
For any $F \in {\cal V}_\phi$,
\[
{\cal S}_\phi{\cal C}_\Omega{\cal S}_\phi^{-1}F(\tau,\sigma) =
\int\int_\Omega F(\tau',\sigma')
\overline{K_{(\tau,\sigma)}(\tau',\sigma')}\,d\tau' d\sigma'.
\]
\label{conc_formula}
\end{proposition}
\noindent{\sc Proof:} Let ${\cal S}_\phi f = F$ and  ${\cal S}_\phi g = G$.
Then
\begin{eqnarray*}
{\cal E}_\Omega (f,g) &=& \int\int_\Omega \langle
f,\rho(\tau,\sigma)\phi\rangle
\langle\rho(\sigma,\tau)\phi,g\rangle\,d\tau d\sigma\\
&=& \int\int_\Omega {\cal S}_\phi f \overline{{\cal S}_\phi g}\,d\tau d\sigma\\
&=& \int\int_\Omega F \overline{G}\,d\tau d\sigma\\
&=& \int\int F|_\Omega \overline{G}\,d\tau d\sigma\\
&=& \int\int {\cal P}_\phi(F|_\Omega) \overline{G}\,d\tau d\sigma.
\end{eqnarray*}
On the other hand, ${\cal E}_\Omega(f,g) = \langle {\cal C}_\Omega
f,g\rangle$. Putting this together easily gives the formula in the
proposition. $\Box$

In view of proposition~\ref{conc_formula}, the spectral properties of the
concentration operator ${\cal C}_\Omega$ are identical to those of the
operator from ${\cal V}_\phi$ to itself defined by
\[
F \mapsto {\cal P}_\phi\chi_\Omega F
\]
where $\chi_\Omega$ is the operator that restricts functions to the set
$\Omega$.

Since the function $\phi$ is fixed in the ensuing discussion, its
use as a subscript will often be supressed in the sequel.

\section{Distribution of the Eigenvalues}

In this section we examine the behavior of the eigenvalues as the
concentration region $\Omega$ gets large. In particular, we show that the
number of eigenvalues that cluster around one grows like the area of the
concentration region. Both the philosophy and proofs in this
section bear a great resemblance to the work of H.J. Landau in
\cite{Landau,Landau2}.

\begin{theorem}
Let $\Omega$ be a bounded measurable set and
\[
n(\lambda ,\Omega_r) = {\rm card}\{k: \lambda_k(\Omega_r) \ge \lambda\}
\]
where $\Omega_r = \{(r\tau ,r\sigma): (\tau,\sigma)\in\Omega\}$. Then
\[
\liminf_{r\to\infty} \frac {n(\lambda ,\Omega_r)}{{\rm area}
(\Omega_r)} = 1
\]
as $r \to \infty$.
\label{thm:asymp_behav}
\end{theorem}

The essential idea behind this theorem, as in the bandlimited case
treated by Landau, is that $\sum \lambda_k$ and $\sum |\lambda_k|^2$
grow at the same asymptotic rate as $r \to \infty$. The following
lemma enables us to work with the operator $\chi_\Omega{\cal P}$ on
$L^2(\Omega)$ as
opposed to the operator ${\cal P}\chi_\Omega$ on ${\cal V}$.

\begin{lemma} The operator $\chi_\Omega{\cal P}$ on
$L^2(\Omega)$ and the operator ${\cal P}\chi_\Omega$ have the
same nonzero eigenvalues with multiplicity.
\label{lem:commute}
\end{lemma}
\noindent{\sc Proof:} Let $F \in {\cal V}$ be a nontrivial
eigenfunction with a nonzero eigenvalue of the operator
${\cal P}\chi_\Omega$:
${\cal P}\chi_\Omega F = \lambda F$. Since both $F$ and
$\lambda$ are nonzero, $\chi_\Omega F \not\equiv 0$. On the
other hand, applying the restriction operator to both sides we
have $\chi_\Omega {\cal P} \chi_\Omega F = \lambda \chi_\Omega
F$. A similar argument shows that every nonzero eigenvalue of
$\chi_\Omega{\cal P}$ is also an eigenvalue of ${\cal
P}\chi_\Omega$.
$\Box$

\begin{lemma}
$\frac {\sum_k |\lambda_k(\Omega_r)|^2}{{\rm area}(\Omega_r)} \to 1$
as $r \to \infty$.
\label{lem:sum_of_squares}
\end{lemma}
\noindent{\sc Proof:} By lemma~\ref{lem:commute}, it suffices
to consider the trace operator $(\chi_{\Omega_r} {\cal P})^2$. The
kernel of this operator on $L^2(\Omega_r)$ is determined by the formula
\[
\begin{array}{rl}
(\chi_{\Omega_r} {\cal P})^2F(\tau,\sigma)&=\\
\int\int\int\int_{\Omega_r\times\Omega_r}& F(\tau'',\sigma'')
\overline{K_{(\tau',\sigma')}(\tau'',\sigma'')}
\overline{K_{(\tau,\sigma)}(\tau',\sigma')}\,d\tau'' d\sigma''
d\tau' d\sigma'.
\end{array}
\]
The kernel is continuous on the bounded set $\Omega_r\times\Omega_r$
so by Mercer's theorem the trace of $(\chi_{\Omega_r} {\cal P})^2$ is
\[
\int\int_{\Omega_r} \int\int_{\Omega_r} |K(\tau,\sigma;\tau',\sigma')|^2\,
d\tau d\sigma d\tau' d\sigma'.
\]
First note that
\begin{eqnarray*}
K(\tau,\sigma;\tau',\sigma') &=& \langle \rho(\tau,\sigma)\phi,
\rho(\tau',\sigma')\phi\rangle\\
&=& \langle \rho(-\tau',-\sigma')\rho(\tau,\sigma)\phi,\phi\rangle\\
&=& \exp\left(\pi\imath (\tau\sigma'-\sigma\tau')\right) \langle
\rho(\tau-\tau',\sigma-\sigma')\phi,\phi\rangle\\
&=& \exp\left(\pi\imath (\tau\sigma'-\sigma\tau')\right)
H(\tau-\tau',\sigma-\sigma').
\end{eqnarray*}
Therefore, it is enough to estimate
\[
\int\int_{\Omega_r}\int\int_{\Omega_r}
|H(\tau-\tau',\sigma-\sigma')|^2\,d\tau d\sigma d\tau' d\sigma'.
\]
By proposition~\ref{spectrogram_identities} the $L^2$ norm of $H$ is
one. The theorem then follows from the following lemma.
$\Box$

\begin{lemma} Let $f \in L^1({\bf R}^n)$ with $f(x) \ge 0$
and $\int f(x) dx = 1$. Let $Q \subset {\bf R^n}$ is a bounded set
of positive measure. Then
\[
r^{-n} \int_{rQ}\int_{rQ} f(x-y) dx\,dy \to |Q|
\]
as $r \to \infty$.
\label{lem:gen_asymp}
\end{lemma}
\noindent{\sc Proof:} Apply the following change of variables:
\[
x = \xi + r\eta \qquad y = r\eta
\]
to get
\begin{equation}
r^{-n} \int_{rQ}\int_{rQ} f(x-y) dx\,dy = \int_Q \int_{r(Q-\eta)}
f(\xi) d\xi d\eta.
\label{eqn:change_of_var}
\end{equation}
Recall that almost every point of a set of
positive measure is a point of density:
\[
\lim_{\epsilon \to 0} \frac {|B_\epsilon(\eta) \cap Q|}{\pi
\epsilon^2} = 1 \qquad \mbox{for almost all $\eta \in Q$.}
\]
Hence for almost all $\eta \in Q$ and arbitrary $R > 0$,
\[
\frac {|r(Q-\eta) \cap B_R(0)|}{\pi R^2} \to 1
\]
as $r \to \infty$. It is then straight-forward to argue that
\[
\int_{r(Q-\eta)\cap B_R(0)} f(\xi) d\xi \to \int_{B_R(0)} f(\xi) d\xi
\]
as $r \to \infty$ for almost all $\eta \in Q$. Since
\[
\int_{r(Q-\eta)\cap B_R(0)} f(\xi) d\xi \le \int_{r(Q-\eta)} f(\xi)
d\xi \le 1
\]
it follows that, for almost all $\eta \in Q$
\[
\int_{B_R(0)} f(\xi) d\xi \le
\liminf_{r\to\infty} \int_{r(Q-\eta)} f(\xi) d\xi \le
\limsup_{r\to\infty} \int_{r(Q-\eta)} f(\xi) d\xi \le 1.
\]
Since $R >0$ was chosen arbitrarily, it follows that
\[
\lim_{r\to\infty} \int_{r(Q-\eta)} f(\xi) d\xi = 1
\]
for almost all $\eta\in Q$. The result then follows from Lebesgue's
dominated convergence theorem.
$\Box$

\noindent{\sc Proof of Theorem~\ref{thm:asymp_behav}:} By Mercer's
theorem, the trace of the operator ${\cal P}\chi$ is
\begin{eqnarray*}
\sum \lambda_k &=& \int\int_{\Omega_r} K(\tau,\sigma,\tau,\sigma)\,d\tau
d\sigma\\
&=& \int\int_{\Omega_r} \langle \rho(\tau,\sigma)\phi,
\rho(\tau,\sigma)\phi\rangle\,d\tau d\sigma\\
&=& {\rm area}(\Omega_r).
\end{eqnarray*}
Choose $r$ so large that $\sum |\lambda_k(\Omega_r)|^2 \ge
(1-\epsilon) {\rm area}(\Omega_r)$. Then
\begin{eqnarray*}
(1-\epsilon) {\rm area}(\Omega_r) &\le& \sum
|\lambda_k(\Omega_r)|^2 \\
&\le& \sum_{k=1}^n \lambda_k + \lambda \sum_{k>n} \lambda_k \\
&\le& (1-\lambda) \sum_{k=1}^n \lambda_k + \lambda
{\rm area}(\Omega_r) \\
&\le& (1-\lambda) n + \lambda {\rm area}(\Omega_r)
\end{eqnarray*}
As a consequence,
\[
n \ge \frac{1-\lambda-\epsilon}{1-\lambda} {\rm area}(\Omega_r).
\]
Since $\epsilon > 0$ was chosen arbitrarily, we have that
\[
\liminf \frac {n(\lambda ,\Omega_r)}{{\rm area}(\Omega_r)} \ge 1.
\]
$\Box$

The following refinement of lemma~\ref{lem:gen_asymp}, in the case
when $Q$ is a bounded domain with $C^1$ boundary and $f$ has
sufficiently strong decay, is useful in order to estimate the size
of the plunge region.

\begin{lemma}
Let $Q$ be a domain with $C^1$ boundary and $f(x)$ as
in lemma~\ref{lem:gen_asymp}. Moreover, assume that
\begin{equation}
|1 - \int_{B_r(0)} f(x)\,dx| \le \frac C{r^p}
\label{eqn:int_decay_cond}
\end{equation}
for some constants $p,C > 0$.
Then the error term
\[
\frac 1{r^n} \int_{rQ}\int_{rQ} f(x-y)\,dx\,dy - |Q|
\]
is ${\cal O}(\frac 1r)$ as $r \to \infty$.
\label{lem:refined_asymp}
\end{lemma}
\noindent{\sc Proof:} Given a $\delta > 0$, let $Q_\delta =
\{x \in Q: {\rm dist}(x,\partial Q) > \delta\}$.
The second iterated integral in equation~\ref{eqn:change_of_var} can
then be expressed as follows:
\begin{equation}
\int_0^{\delta_0} \int_{\partial Q_\rho} I_r(\eta)\,d{\cal
H}(\eta)\,d\rho
\label{eqn:iterated_int}
\end{equation}
where
\[
I_r(\eta) = \int_{r(Q-\eta)}f(\xi)\,d\xi,
\]
$d{\cal H}$ is Hausdorff $n-1$ dimensional measure, and $\delta_0$
is any positive number larger than the diameter of $Q$. Using
equation~\ref{eqn:int_decay_cond} to estimate $I_r(\eta)$ yields
\begin{eqnarray}
\left| \int_{1/r}^{\delta_0} \int_{\partial Q_\rho} I_r(\eta)\,d{\cal
H}(\eta)\,d\rho - |Q_{1/r}|\right|
& \le & \int_{1/r}^{\delta_0} \int_{\partial Q_\rho} |I_r(\eta)-1|
\,d{\cal H}(\eta)\,d\rho\nonumber\\
& \le &
C \int_{1/r}^{\delta_0} \int_{\partial Q_\rho} \frac 1{(r\rho )^p}
\,d{\cal H}(\eta)\,d\rho.\label{eqn:piec1}
\end{eqnarray}
Using the fact that the Hausdorff $n-1$ measure of the sets
$\partial Q_\rho$ are uniformly bounded, one easily has an estimate of the
required type for the first expression in
inequality~\ref{eqn:piec1}. Moreover, since $0 \le I_r(\eta) \le 1$,
\begin{equation}
\begin{array}{rl}
|\int_0^{1/r} \int_{\partial Q_\rho} I_r(\eta)\,d{\cal
H}(\eta)\,d\rho - {\rm meas}(Q\backslash Q_{1/r})|
& = \\
\int_0^{1/r} \int_{\partial Q_\rho} &
(1 - I_r(\eta))\,d{\cal
H}(\eta)\,d\rho \le\\
\int_0^{1/r} \int_{\partial Q_\rho}&\,d{\cal
H}(\eta)\,d\rho \\
\le &
{\rm meas}(Q\backslash Q_{1/r}) \le C/r
\end{array}
\label{eqn:piece2}
\end{equation}
Putting equations~\ref{eqn:piec1} and \ref{eqn:piece2}
together via the triangle inequality yields the sought after result.
$\Box$

\begin{theorem}
Let $\Omega$ be a bounded domain with $C^1$ boundary and
\[
n(\lambda,\mu,\Omega_r) = {\rm card}\{k:\lambda_k(\Omega_r) \in
[\lambda,\mu]\}
\]
where $0 < \lambda < \mu < 1$. Moreover, assume that the function
$|{\cal S}_\phi\phi(\tau,\sigma)|^2$ satisfies
condition~\ref{eqn:int_decay_cond}. Then, there is a constant $C > 0$
for which
\[
n(\lambda,\mu,\Omega_r) \le C r.
\]
\label{thm:plunge}
\end{theorem}
\noindent{\sc Proof:} Mercer's theorem applied to the positive
operator $\chi_{\Omega_r}{\cal P} - (\chi_{\Omega_r}{\cal P})^2$
yields that
\[
\sum \lambda_k(\Omega_r) - \lambda_k(\Omega_r)^2 = |\Omega_r| -
\int\int_{\Omega_r}\int\int_{\Omega_r}
|H(\tau-\tau',\sigma-\sigma')|^2\,d\tau\, d\tau'\, d\sigma\,
d\sigma'.
\]
Applying lemma~\ref{lem:refined_asymp} yields
\[
0 < \sum \lambda_k(\Omega_r) - \lambda_k(\Omega_r)^2 \le Cr.
\]
Therefore, if $\epsilon \in (0,1/2)$ then
\[
n(\epsilon,1 - \epsilon,\Omega_r) (\epsilon - \epsilon^2) \le
\sum \lambda_k(\Omega_r) - \lambda_k(\Omega_r)^2 \le Cr.
\]
The result follows.
$\Box$

\section{Regularity of the Eigenfunctions}

In this section we will be concerned with the regularity of the
eigenfunctions of the concentration operator ${\cal C}_\Omega$. In
order to do this it will be necessary to work with the integral
kernel of this operator.

\begin{proposition}
The kernel $k_\Omega(x,y)$ of the operator ${\cal C}_\Omega$ is
\[
k_\Omega(x,y) = \int\int_\Omega \rho(\tau,\sigma)\phi(x)
\overline{\rho(\tau,\sigma)\phi(y)}\,d\tau d\sigma.
\]
\label{prop:kernel_formula}
\end{proposition}
{\sc Proof:} Using the definition of the concentration operator,
one has
\begin{eqnarray*}
\langle {\cal C}_\Omega f,g\rangle &=& {\cal E}_\Omega(f,g) \\
&=& \int\int_\Omega \langle f,\rho(\tau,\sigma)\phi\rangle\,
\langle \rho(\tau,\sigma)\phi,g\rangle\,d\tau d\sigma \\
&=& \int\int f(y) \left(\int\int_\Omega \rho(\tau,\sigma)\phi(x)
\overline{\rho(\tau,\sigma)\phi(y)}\,d\tau d\sigma\right)
\overline{g(x)}\,dx dy.
\end{eqnarray*}
The result follows.
$\Box$

Suppose that $\gamma(s)$ is a positive decreasing function defined
on $[0,\infty)$ satisfying the following conditions:
\begin{itemize}
\item $\gamma(0) = 0$ and $\lim_{s\to\infty} \gamma(s) = 0$,
\item there is an $\epsilon_0 > 0$ such that
$\int |\gamma(s)|^p\,ds < \infty$ for all $p \in (2 - \epsilon_0, 2]$
, and
\item for any $s_0 > 0$, $\gamma(s - s_0) = {\cal
O}(\gamma(s)^{1-\epsilon})$ for all sufficiently small
$\epsilon > 0$.
\end{itemize}
Examples of functions $\gamma(s)$ that satisfy such conditions
include
\begin{itemize}
\item $\gamma(s) = (1 + s^2)^{-q/2}$ where $q \in (1,\infty)$, and
\item $\gamma(s) = \exp(- \kappa |s|^q)$ where $\kappa,q > 0$.
\end{itemize}

\begin{proposition}
If $|\phi(t)| \le \gamma(|x|)$ for all $t$ then
any eigenfunction with nonzero eigenvalue $\psi$ of the
concentration operator ${\cal C}_\Omega$ satisfies an estimate of
the form
\[
|\psi(t)| \le C \gamma(|t|)^{1-\epsilon}
\]
for all $t$.
\end{proposition}
{\sc Proof:} Suppose that ${\cal C}_\Omega\psi = \lambda \psi$ with
$\lambda > 0$. The first step is to estimate the kernel,
$k_\Omega(x,y)$:
\begin{eqnarray*}
|k_\Omega(x,y)| &=& \left|\int\int_\Omega e^{2\pi\imath \sigma(x-y)}
\phi(x+\tau)\overline{\phi(y+\tau)}\,d\tau d\sigma\right|\\
&\le & C \gamma(\left||x| - R\right|) \gamma(\left||y| - R\right|)\\
&\le& C \gamma(|x|)^{1-\epsilon} \gamma(|y|)^{1-\epsilon}
\end{eqnarray*}
where $R = \sup \{|(\tau,\sigma)|:(\tau,\sigma)\in\Omega\}$.
Consequently, using the Cauchy-Schwarz inequality, we have the
estimate:
\begin{eqnarray*}
\lambda |\psi(x)| &\le & C \gamma(|x|)^{1-\epsilon} \int
\gamma(|y|)^{1-\epsilon} \psi(y)\,dy\\
&\le & C \gamma(|x|)^{1-\epsilon}.
\end{eqnarray*}
$\Box$

Let $\tilde\Omega = \{(\tau,\sigma):(-\sigma,\tau)\in\Omega\}$.
The second item in proposition~\ref{spectrogram_identities} implies that:
\begin{eqnarray*}
{\cal E}_\Omega(f,g) &=& \int\int_\Omega {\cal S}_\phi f
\overline{{\cal S}_\phi g} \\
&=& \int\int_{\tilde\Omega} {\cal S}_{\hat \phi} \hat f
\overline{{\cal S}_{\hat \phi} \hat g}\\
&=& {\cal E}_{\tilde\Omega}(\hat f,\hat g).
\end{eqnarray*}
In particular,
\[
{\cal E}_{\tilde\Omega}(\hat \psi_k,\hat \psi_l) = \delta_{kl}
\]
and $\hat \psi_1,\cdots$ are the eigenfunction of the operator
${\cal C}_{\tilde\Omega}$. Putting this together with the preceding
proposition yields the following theorem.

\begin{theorem} If
\begin{eqnarray*}
|\phi(t)| & \le & \gamma_1(|t|) {\  and}\\
|\hat\phi(\sigma)| &\le& \gamma_2(|\sigma|)
\end{eqnarray*}
where both $\gamma_1$ and $\gamma_2$ satisfy the conditions
described above, then
\begin{eqnarray*}
|\psi(t)| & \le & C \gamma_1(|t|)^{1-\epsilon} {\ and}\\
|\hat\psi(\sigma)| & \le & C \gamma_2(|\sigma|)^{1-\epsilon}.
\end{eqnarray*}
\label{thm:reg}
\end{theorem}

Various special cases merit mention.
\begin{itemize}

\item Let
\[
\phi(t) = \left\{
\begin{array}{ll}
1 - |t| & \mbox{if $|t| \le 1$}\\
0 & \mbox{otherwise}.
\end{array}
\right\}
\]
Then $k_\Omega(x,y) = 0$ whenever $|x-y| \ge 2$. Hence, any
eigenfunction, with nonzero eigenvalue, must have support contained
in the interval $[-3,3]$. Moreover, the Fourier transform of $\phi$
is the Fejer kernel which is bounded by $\gamma(s) = (1 + s^2)^{-1}$. As a
consequence, if
Hence, if $\psi$ is any eigenfunction with nonzero eigenvalue,
\[
\begin{array}{ccl}
\mbox{support($\psi$)} &\subset& [-3,3]\\
\hat \psi(\sigma) &=& {\cal O}(\frac 1{1+|\sigma|^{2-\epsilon}}).
\end{array}
\]
In particular, $\psi'$ is square integrable.

\item If $\phi$ is Schwartz class, one can verify directly using
proposition~\ref{prop:kernel_formula} that the kernel
$k_\Omega$ is also Schwartz class. A straightforward interchange of the
differentiation and integration symbol then
shows that $\psi$ is also Schwartz class:
\[
\begin{array}{ccl}
|\lambda t^\alpha D_t^\beta \psi(t)| &=& |\int t^\alpha D_t^\beta
k_\Omega(t,s) \psi(s) ds|\\
&\le& \|t^\alpha D_t^\beta k_\Omega(t,s)\|\ \|\psi\|.
\end{array}
\]
Since the kernel is Schwartz class, the last expression is bounded
independently of $t$.  Thus, if $\phi$ is Schwartz class,
so are any eigenfunctions associated to nonzero eigenvalues.

\item If $\phi(t) = \exp(- C t^2)$, theorem~\ref{thm:reg} implies that the
eigenfunctions of the concentration operator with nonzero
eigenvalue are
${\cal O}(\exp(- C(1-\epsilon) t^2))$. Moreover, since
$\hat\phi(\sigma) = \sqrt{\pi/C} \exp(-\pi/C \sigma^2)$, the Fourier
transforms of these eigenfunctions are ${\cal O}(\exp(-\frac
{\pi}{C} (1-\epsilon) \sigma^2)$.  In particular, these functions are
entire by the Paley-Wiener theorem ~\cite{Katznelson}.  Furthermore, by
differentiating under the integral sign (as above, and in the formula for the
kernel
in proposition~\ref{prop:kernel_formula}), one can show that all derivatives of
these
eigenfunctions satisfy the same decay estimates.

\end{itemize}

Theorem~\ref{thm:reg} is useful in that
it provides quantitative bounds on the decay (and eventually
regularity) of the eigenfunctions for a wide class of $\phi$.  Among the works
in
the literature that bear on the problems discussed here, we remark that
Pietsch~\cite{Pietsch} provides estimates on the decay of the eigenvalues for
given
decay and regularity of the integral kernel.  In particular, his results imply
that
when $\phi$ is in Schwartz space, the sequence of eigenvalues decays faster
than the
inverse of any polynomial in the index n.  Janssen~\cite{Janssen1} obtains the
same estimate
in the Schwartz class case.

For $\phi(t) = \exp(- C t^2)$, the work of Daubechies~\cite{Daubechies} is
relevant.  She solves the eigenvalue problem explicitly for circular (or
elliptical) domains, and obtains the Hermite functions as eigenfunctions.  She
also
obtains a formula for the eigenvalues, which decay exponentially in the
index n.  Our approach shows that for any bounded measurable domain, the
eigenfunctions are analytic and have quadratic exponential decay on the real
line.  Since the
kernel is Schwartz class, ~\cite{Pietsch} implies that the eigenvalues decay
faster than the
inverse of any polynomial in the index n.  In fact, a stronger statement can be
obtained
using a result of Janssen~\cite{Janssen2}.  Note that the integral kernel is
\[
\begin{array}{ccl}
k_\Omega(x,y) = \int\int_\Omega exp(2\pi i \sigma (x-y)) exp(- C
((x+\tau)^2 + (y+\tau)^2))\,d\tau d\sigma.
\end{array}
\]

By expanding the exponentials in power series and integrating term by term, it
is
easy to show that this kernel is analytic in each variable, and satisfies the
bound
${\cal O}(\exp(- A Re[x]^2 + B Im[x]^2 - C Re[y]^2 + D Im[y]^2))$, for some
$A, B, C, D > 0$.  Since the kernel is positive definite as well,
it follows from theorem A.1 in ~\cite{Janssen2} that the eigenfunctions
belonging to
nonzero eigenvalue are analytic and satisfy ${\cal O}(\exp(- A Re[x]^2
+ B Im[x]^2))$, some  $A, B > 0$.  At
the same time, the eigenvalues are in ${\cal O}(\exp(- n \alpha))$, some
$\alpha > 0$.


\begin{thebibliography}{99}

\bibitem{Daubechies}
I. Daubechies:
Time-Frequency Localization Operators: A Geometric Phase Space
Approach I,
{\em IEEE Trans. on Information Theory},
{\bf 34} (1988), pp. 605-612.

\bibitem{Daubechies and Paul}
I. Daubechies and T. Paul:
Time-Frequency Localization Operators: A Geometric Phase Space
Approach II,
{\em Inverse Problems},
{\bf 4} (1988), pp. 661-680.

\bibitem{Folland}
G. B. Folland:
{\em Harmonic Analysis in Phase Space},
Princeton University Press,
Princeton, NJ, 1989.

\bibitem{Janssen1}
A. J. E. M. Janssen:
On Hermitian Positive Definite Functions of Two Variables,
unpublished.

\bibitem{Janssen2}
A. J. E. M. Janssen:
Positivity Properties of Phase-Plane Distribution Functions,
{\em J. Math. Phys.}, {\bf 25} (1984), pp. 2240-2252.

\bibitem{Katznelson}
Y. Katznelson:
{\em An Introduction to Harmonic Analysis},
John Wiley and Sons,
New York, NY, 1968.

\bibitem{Landau}
H.J. Landau:
Necessary Density Conditions for Sampling and Interpolation of
Certain Entire Functions,
{\em Acta Mathematica}, {\bf 117} (1967), pp. 37-52.

\bibitem{Landau2}
H.J. Landau:
On Szego's Eigenvalue Distribution Theorem and Non-Hermitean
Kernels,
{\em J. Analyse Math.}, {\bf 28} (1975), 335-357.

\bibitem{Pietsch}
A. Pietsch:
Eigenvalues of Integral Operators II,
{\em Mathematische Annalen},
{\bf 262} (1983), 343-376.

\end{thebibliography}
\end{document}